\documentclass[a4paper,titlepage]{article}
\pdfoutput=1

\usepackage[utf8]{inputenc}
\usepackage[T1]{fontenc}
\usepackage[UKenglish]{babel}
\usepackage{ifthen}
\usepackage{tabulary}
\usepackage{chngcntr}

\usepackage{amssymb}
\usepackage{amsmath}
\usepackage{amsthm}
\usepackage{array}
\usepackage{eqnarray}
\usepackage{mathtools}
\usepackage[ruled, vlined, linesnumbered]{algorithm2e}
\let\chapter\undefined			
\usepackage{listings}
\lstdefinelanguage{Mosel}{
  morekeywords={array, as, boolean, break, case, count, counter, declarations, do, dynamic, elif, else, evaluation, exists, false, forall, forward, from, function, if, imports, include, initialisations, initializations, integer, inter, is_binary, is_continuous, is_free, is_integer, is_partint, is_semcont, is_semint, is_sos1, is_sos2, linctr, list, model, mpvar, next, of, options, package, parameters, procedure, public, prod, range, real, record, repeat, requirements, set, string, sum, then, to, true, union, until, uses, version, while, while, end-do,
end-if, end-initialisations, end-initializations, end-model, end-package, end-declarations, end-requirements, end-parameters, end-record, end-function, end-procedure},
  morekeywords=[2]{and, xor, or, not, mod, div, max, min},
  morekeywords=[3]{=,:=,+=,-=,>,>=,<,<=,<>,+,-,*,/,in},
  morekeywords=[4]{finalize, finalise, create, getsize, writeln, write, getparam, setparam, exit, wait, sleep, minimize, maximize, minimise, maximise, text, minlist, maxlist},
  morekeywords=[5]{SQLconnect, SQLexecute, SQLdisconnect, time, date, datetime, getdate, setcallback, XPRS_CB_OPTNODE, XPRS_INT_SOL, XPRS_DUAL},
  sensitive=false,
  alsoletter={+-*/<>=:},
  morecomment=[l]!,
  morecomment=[s]{(!}{!)},
  morestring=[b]",
  morestring=[b]',
  tabsize=4
}
\lstnewenvironment{Mosel}[1]
{\listingfile{#1}\nopagebreak\lstset{language=Mosel}}
{}

\newcommand{\mosel}[1]{{\lstset{language=Mosel}\lstinline$#1$}}

\newcommand{\listingfile}[1]{\ifthenelse{\equal{#1}{}}{}%
{%
\hfill\hskip3ex\vskip-\baselineskip\vskip\topsep\hbox{\ttfamily\hfill\footnotesize #1}\vskip-\topsep
}}
\usepackage{multirow}

\usepackage{pgf}
\usepackage{tikz}
\usetikzlibrary{shapes}
\usetikzlibrary{decorations.pathmorphing}
\usepackage{graphicx}
\usepackage{floatflt}
\usepackage{float}
\usepackage[normalsize, hang]{subfigure}
\usepackage[normal]{caption}

\usepackage{times}
\usepackage[T1]{fontenc}
\usepackage{microtype}
\usepackage{eurosym}

\usepackage{longtable}
\usepackage{setspace}             
\numberwithin{equation}{section}  
\makeatletter
\renewcommand{\fps@figure}{htb}   
\renewcommand{\fps@table}{htb}    
\makeatother
\usepackage{mdwlist}              
\usepackage{fancybox}             %

\newcommand{\shortcaption}[2][\empty]{{%
	\renewcommand{\figurename}{Fig.}%
	\ifthenelse{\equal{#1}{\empty}}%
		{\caption{#2}}%
		{\caption[#1]{#2}}%
}}

\newcommand{\set}[1]{\mathbb{#1}}

\def\discrange#1#2{{%
	\toks255={#1}
	\ifthenelse{\equal{\the\toks255}{1}}%
		{\lbrack#2\rbrack}%
		{\lbrack#1..#2\rbrack}%
}}

\renewcommand{\quote}[1]{``#1''}

\makeatletter
\newcommand{\newvar}[3]{%
	\expandafter\def\csname #1\endcsname{\expandafter\@ifnextchar\expandafter(\csname #1Params\endcsname{\error{This variable should have parameters!}}}%
	\ifcase #2
		\expandafter\def\csname #1\endcsname{#3\@ifnextchar({\error{This variable should not have parameters!}}\relax}%
		\expandafter\def\csname #1b\endcsname{#3}%
	\or	\expandafter\def\csname #1Params\endcsname(##1){#3}%
		\expandafter\def\csname #1b\endcsname##1{#3}%
	\or	\expandafter\def\csname #1Params\endcsname(##1,##2){#3}%
		\expandafter\def\csname #1b\endcsname##1##2{#3}%
	\or	\expandafter\def\csname #1Params\endcsname(##1,##2,##3){#3}%
		\expandafter\def\csname #1b\endcsname##1##2##3{#3}%
	\or	\expandafter\def\csname #1Params\endcsname(##1,##2,##3,##4){#3}%
		\expandafter\def\csname #1b\endcsname##1##2##3##4{#3}%
	\or	\expandafter\def\csname #1Params\endcsname(##1,##2,##3,##4,##5){#3}%
		\expandafter\def\csname #1b\endcsname##1##2##3##4##5{#3}%
	\or	\expandafter\def\csname #1Params\endcsname(##1,##2,##3,##4,##5,##6){#3}%
		\expandafter\def\csname #1b\endcsname##1##2##3##4##5##6{#3}%
	\or	\expandafter\def\csname #1Params\endcsname(##1,##2,##3,##4,##5,##6,##7){#3}%
		\expandafter\def\csname #1b\endcsname##1##2##3##4##5##6##7{#3}%
	\or	\expandafter\def\csname #1Params\endcsname(##1,##2,##3,##4,##5,##6,##7,##8){#3}%
		\expandafter\def\csname #1b\endcsname##1##2##3##4##5##6##7##8{#3}%
	\or	\expandafter\def\csname #1Params\endcsname(##1,##2,##3,##4,##5,##6,##7,##8,##9){#3}%
		\expandafter\def\csname #1b\endcsname##1##2##3##4##5##6##7##8##9{#3}%
	\fi
}
\makeatother

\usepackage[pdfpagelabels]{hyperref}

\newcommand{\Graphics}{.}

\newvar{periodLength}{0}{\text{L}}
\newvar{minProd}{1}{\underline{\text{P}}_{#1}}
\newvar{maxProd}{1}{\overline{\text{P}}_{#1}}
\newvar{minUptime}{1}{\text{UT}_{#1}}
\newvar{minDowntime}{1}{\text{DT}_{#1}}
\newvar{initialOffline}{1}{\text{IDT}_{#1}}
\newvar{initialOnline}{1}{\text{IUT}_{#1}}
\newvar{rampup}{1}{\text{RU}_{#1}}
\newvar{rampdown}{1}{\text{RD}_{#1}}
\newvar{startupRampup}{1}{\text{SU}_{#1}}
\newvar{shutdownRampdown}{1}{\text{SD}_{#1}}
\newvar{demand}{1}{\text{D}^{#1}}
\newvar{reserve}{1}{\text{R}^{#1}}
\newvar{storage}{1}{\text{S}}
\newvar{initialStorage}{1}{\text{SI}_{#1}}
\newvar{finalStorage}{1}{\text{SF}_{#1}}
\newvar{storageCapacity}{1}{\text{SC}_{#1}}
\newvar{storageEfficiency}{1}{\text{SE}_{#1}}
\newvar{storageInflow}{1}{\text{SIF}_{#1}}
\newvar{StartupCost}{2}{\text{CU}_{#1}^{#2}}
\newvar{ShutdownCost}{1}{\text{CD}_{#1}}
\newvar{underProductionPenalty}{0}{\text{UPP}}
\newvar{underReservePenalty}{0}{\text{URP}}
\newvar{overProductionPenalty}{0}{\text{OPP}}
\newvar{fuelType}{1}{\text{F}_{#1}}
\newvar{fuelCost}{2}{\text{FC}_{#1}^{#2}}
\newvar{varFuel}{1}{\text{FA}_{#1}}
\newvar{fixedFuel}{1}{\text{FB}_{#1}}
\newvar{varCost}{1}{\text{PA}_{#1}}
\newvar{fixedCost}{1}{\text{PB}_{#1}}

\newvar{onOff}{2}{v_{#1}^{#2}}
\newvar{prod}{2}{p_{#1}^{#2}}
\newvar{price}{1}{\text{price}^{#1}}
\newvar{consumption}{2}{c_{#1}^{#2}}
\newvar{storage}{2}{s_{#1}^{#2}}
\newvar{maxPossProd}{2}{\overline{p}_{#1}^{#2}}
\newvar{underProduction}{1}{P_-^{#1}}
\newvar{overProduction}{1}{P_+^{#1}}
\newvar{underReserve}{1}{R_-^{#1}}
\newvar{prodCost}{2}{\text{cp}_{#1}^{#2}}
\newvar{startupCost}{2}{\text{cu}_{#1}^{#2}}
\newvar{shutdownCost}{2}{\text{cd}_{#1}^{#2}}

\newcommand{\forallK}{\forall\ k \in K}
\newcommand{\forallJK}{\forall\ j \in J,\, k \in K}
\newcommand{\forallJ}{\forall\ j \in J}

\hypersetup{
    pdftitle={Implementing a Unit Commitment formulation},
    pdfkeywords={Unit Commitment} {Power production} {key3},
    pdfauthor={Matthias Silbernagl - Technische Universit\"at M\"unchen},
    colorlinks=false,       
    linkcolor=red,          
    citecolor=green,        
    filecolor=magenta,      
    urlcolor=cyan,          
    plainpages=false
}

\begin{document}


\title{Implementing a Unit Commitment Power Market Model in FICO Xpress-Mosel\thanks{Preprint; Published version available at \url{http://community.fico.com/docs/DOC-1365}}}
\author{Ren{\'e} Brandenberg \\ Technische Universit\"at M\"unchen \and Matthias Silbernagl \\Technische Universit\"at M\"unchen}
\date{}
\renewcommand{\thepage}{Title}
\maketitle

\renewcommand{\thepage}{Abstract}
\begin{abstract}
Starting from a basic Unit Commitment problem as published in
\cite{carrion2006computationally}, we develop a fully functional, running
implementation for Xpress Mosel, useable in power market modeling.
The constraints of the model are discussed in details and solutions to left open 
implementation problems are presented, guidelines on how to handle data input 
(from both, databases and Excel), as well as data output (to Excel) are given.
\end{abstract}

	
\pagenumbering{Roman}
\tableofcontents
\clearpage

	
\pagenumbering{arabic}
\setcounter{page}{1}

\section{Introduction}

The deregulation of the energy market in the last two decades attracted the researchers' attention, as electricity producers constantly needed to optimize their production to stay competitive.
The raising volatility in energy production, mostly due to the renewable resources, even increases the potential for optimization and therefore the interest in power market models.

We consider a Unit Commitment (UC) model based on the mixed integer linear program (MILP) presented in \cite{carrion2006computationally}, which assumes perfect competition, and as such approximates the typical oligopolistic power market in deregulated countries. Its main advantage over possibly more sophisticated game theoretical or stochastic approaches is its capability to handle real-world problem sizes with an accurate physical model of the power units, combined with its relatively good comprehensibility. 

To make the model in \cite{carrion2006computationally} suitable for real-life application, we need to modify it in a few key aspects. In particular, we focus on
\begin{itemize}
  \item the reduction of the prohibitively large number of constraints needed to model the start-up costs,
  \item the incorporation of storage units,
  \item the analysis of causes of infeasibility and how to avoid them, and
  \item the incorporation of daily fuel prices, as fuel costs make up the better part of the production costs.
\end{itemize}

After a short theoretical presentation of the model, we put our main focus on a discussion of the constraints and a thorough guide on their implementation.

Finally, we show how to read the needed input data (which may be considerably huge, if long time period real world data is used) from either an SQL database or an Excel spreadsheet. The output of the results is written to an Excel spreadsheet, for further processing and analysis.

We use the Xpress-Mosel modeling language \cite{mosel}, which is part of the Xpress Optimization Suite \cite{xpress}.

\clearpage

\section{The Unit Commitment Problem}

In the Unit Commitment problem one has to satisfy a certain energy demand over a number of periods, using a set of given power units, while minimizing the total incurred costs.

A typical Unit Commitment problem can roughly be split into three parts:
\begin{enumerate}
  \item \textbf{The modeling of the physical constraints of each unit (subsection \ref{sec:PhysicalConstraints})}\\
    In this part the technical limitations of each individual unit are
    considered. We include
    \begin{itemize}
      \item the feasible range for the production level,
      \item the feasible range for the change rate of the production level,
      \item the minimal up- and downtime and
      \item the storage constraints.
    \end{itemize}
  \item \textbf{The modeling of the power grid (subsection \ref{sec:DemandAndReserve})}\\
    A model of the undelying power grid may be considered, too, which we exclude in this paper. Hence, the overall production is just connected to the demand, a typical assumption in many models, usually described by the keyword \quote{big copper plate}.  
  \item \textbf{The modeling of the costs (subsection \ref{sec:CostConstraints})}\\
    Besides the obvious fuel costs, we include start up and shutdown costs and some penalties.
\end{enumerate}

\noindent In the remaining of this section we discuss the components of our
model in the above listed order (subsections \ref{sec:PhysicalConstraints} -
\ref{sec:CostConstraints}), give a short note on how to handle infeasibilities
\ref{sec:DealingWithInfeasibilities}, and include summaries of the parameters
\ref{sec:Parameters} and the model \ref{sec:Model}. The latter may be used as a convenient reference in the following sections, since it is cross-referenced with the discussion and implementation over there. 

\subsection{Basic notation}

We denote the set of time periods as $K = \{1, \ldots, T\}$, and the set of units as $J$. The remaining parameters will be introduced ``as we go'', and summarized in subsection \ref{sec:Parameters}.
We use the shortcut notations $\discrange{a}{b}$ to denote the discrete set $[a, b] \cap \set{Z}$, and $\discrange{1}{b}$ to denote $\discrange{1\relax}{b}$.


\subsection{Physical Constraints of the Units}
\label{sec:PhysicalConstraints}

\subsubsection{Unit Variables (constraints \ref{eq:UnitVariablesStart} and \ref{eq:UnitVariablesEnd} in subsection \ref{sec:Model})}
\label{sec:UnitVariables}

In any period $k \in K$, the units $j \in J$ are modeled by their operational state (on/off) $\onOff(j, k)$, their current production $\prod(j, k)$ and their maximal possible production $\maxPossProd(j,k)$
\begin{alignat*}{2}
	& \forallJK: \qquad & \onOff(j,k) &\in \{0, 1\}\\
	&                   & \prod(j,k), \maxPossProd(j,k) &\in \set{R}^+,
\end{alignat*}
where $\onOff(j,k) = 1$ iff unit $j$ is operational in period $k$ and may produce power.

\clearpage\noindent The maximal possible production $\maxPossProd(j,k)$ is used to measure the
\textit{spinning reserve}: In case of a power outage, the power grid must be
stabilized by ramping up the currently operational units. The available
additional production capacity $\maxPossProd(j, k) - \prod(j, k)$ is called
\textit{spinning reserve} of unit $j$ and a general grid constraint regulates the minimal need of spinning reserve
over all units (see subsection \ref{sec:DemandAndReserve} below).


\subsubsection{Minimal Up- and Downtime (constraints \ref{eq:MinimalUpAndDowntimeStart} to \ref{eq:MinimalUpAndDowntimeEnd} in subsection \ref{sec:Model})}
\label{sec:MinimalUpAndDowntime}

Units cannot be started up and shut down arbitrarily. For example, for a coal unit, after being shut down, an appropriate cool-off time has to be kept.

To model this, minimal up- and downtime parameters $\minUptime(j)$ and $\minDowntime(j)$ for each unit are given, denoting the minimal time a unit has to stay operational after a start up, and the minimal cool-off time after a shutdown, respectively. Since we do not know when a unit has last been started up or shut down prior to the modeled time range, we further expect to be given the two parameters $\initialOnline(j)$ and $\initialOffline(j)$, which tell us for how long a unit needs to be operational or shut down initially. Given these parameters, we model the minimal uptime as
\[\forallJ,\, k \in \discrange{1}{\initialOnline(j)}: \qquad \onOff(j,k) = 1\]
\vskip-0.7em\noindent for the first periods, and as\vskip-1.9em
\[\forallJ,\, k \in \discrange{(\initialOnline(j) + 2)}{T},\, i \in \discrange{1}{\minUptime(j) - 1} \cap \discrange{1}{T - k}:\qquad\onOff(j,k+i) \geq \onOff(j,k) - \onOff(j,k-1)\]
\noindent for the remaining periods. Consider the right-hand side $\onOff(j, k) - \onOff(j, k-1)$ of the latter constraint: it always lies in $\{-1, 0, 1\}$, and is equal to $1$ if and only if unit $j$ starts up in period $k$. Thus, in the start-up case, the variables $\onOff(j,k+i)$ on the left-hand side are forced to $1$, while the constraint is always fulfilled otherwise. Similarly, the minimal downtime is modeled as
\[\forallJ,\, k \in \discrange{1}{\initialOffline(j)}: \qquad \onOff(j,k) = 0\]
\[\forallJ,\, k \in \discrange{(\initialOffline(j)\!+\!2)}{T},\, i \in \discrange{1}{\minDowntime(j) - 1} \cap \discrange{1}{T - k}:\qquad \onOff(j,k+i) \leq 1 - (\onOff(j,k) - \onOff(j,k-1))\]
\noindent Note that we expect $\initialOffline(j), \initialOnline(j) \in
\discrange{0}{T}$ and $\minDowntime(j), \minUptime(j) \in \discrange{1}{T}$.


\subsubsection{Minimal and Maximal Production (constraint \ref{eq:MinimalAndMaximalProduction} in subsection \ref{sec:Model})}
\label{sec:MinimalAndMaximalProduction}

If a unit $j$ is not operational, it may not produce; otherwise both the actual and maximal possible production level have to lie in a specific interval, defined by the two parameters $\minProd(j)$ (minimal production) and $\maxProd(j)$ (maximal production). All this can be expressed by
\[\forallJK: \qquad \minProd(j) \onOff(j,k) \ \leq\  \prod(j,k) \ \leq\  \maxPossProd(j,k) \ \leq\  \maxProd(j) \onOff(j,k)\]
Note that storage units may have a negative minimal production
$\minProd(j)$. However, $\prod(j,k) \ge 0$ only captures power
production, while power consumption is modeled by $\consumption(j,k)$ (see
section \ref{sec:Storage} below).


\subsubsection{Ramping (constraints \ref{eq:RampingStart} to \ref{eq:RampingEnd} in subsection \ref{sec:Model})}
\label{sec:Ramping}

Of course the production level can not change arbitrarily either. The parameters
$\rampup(j)$ / $\rampdown(j)$ denote the maximal speed when increasing /
decreasing the production level of an operating unit. For example, a unit with
$\rampup(j) = 50\text{MW/h}$ would be able to increase it's production level
from $120\text{MW}$ to $170\text{MW}$ in one hour, but not to
$171\text{MW}$. This may simply be modeled as
\begin{alignat*}{2}
	& \forallJ,\, k \in \discrange{2}{T}:\qquad & \maxPossProd(j,k) & \leq \prod(j,k-1) + L \cdot \rampup(j)\\
	&                                           & \prod(j,k) & \geq \prod(j,k-1) - L \cdot \rampdown(j) \onOff(j,k)
\end{alignat*}

\noindent where $L$ is a parameter denoting the period length. However, at start up and shutdown units are typically able to change their production levels faster. This higher ramping speed is denoted by the two parameters $\startupRampup(j)$ (maximal production level at start up) and $\shutdownRampdown(j)$ (maximal production level before shutdown).

Thus for up ramping we get the constraint
\begin{alignat*}{2}
	& \forallJ,\, k \in \discrange{2}{T}:\qquad & \maxPossProd(j,k) & \leq \prod(j,k-1) + L \cdot \rampup(j) \onOff(j,k-1) + \startupRampup(j) (1 - \onOff(j,k-1))\\
	&&&- \min\{\startupRampup(j), \max\{\minProd(j), 0\} + L \cdot \rampup(j)\} \cdot (1 - \onOff(j,k))
\end{alignat*}
The first three terms on the right-hand side ascertain that the production may only increase by $L \cdot \rampup(j)$ if the unit is already running, or by $\startupRampup(j)$ if the unit is starting.

The last term does not change the right hand side in the case $\onOff(j,k) = 1$, but tightens the constraint in the case $\onOff(j,k) = 0$. Consider the constraint without the last term:
\begin{alignat*}{4}
\textrm{Case }\onOff(j,k) = 0,\, \onOff(j,k-1) = 0:\quad& \underbrace{\maxPossProd(j,k)}_{= 0} &&\leq \underbrace{\prod(j,k-1)}_{= 0} + &&\underbrace{L \cdot \rampup(j) \onOff(j,k-1)}_{= 0} + \startupRampup(j) \underbrace{(1 - \onOff(j,k-1))}_{=1} &&= \startupRampup(j)\\
\textrm{Case }\onOff(j,k) = 0,\, \onOff(j,k-1) = 1:\quad&
\underbrace{\maxPossProd(j,k)}_{= 0} &&\leq \underbrace{\prod(j,k-1)} + \ &&L \cdot \rampup(j) \underbrace{\onOff(j,k-1)}_{= 1} + \underbrace{\startupRampup(j) (1 - \onOff(j,k-1))}_{=0} &&= \prod(j,k-1) + L \cdot \rampup(j)
\end{alignat*}
Thus, in case of $\onOff(j,k) = 0$ we could tighten the constraint by
subtracting the term $\min\{\startupRampup(j), \prod(j,k-1) + L \cdot
\rampup(j)\}$ on the right hand side. However, since we need to multiply the
term by $(1 - \onOff(j,k))$ to leave the case $\onOff(j,k) = 1$ unchanged, to
avoid non-linearities, it should not contain a variable. Hence we have to replace $\prod(j,k-1)$ by its best lower bound $\max\{\minProd(j), 0\}$. This leads exactly to the last term in the rampup constraint. This term is not necessary for a correct model, but tightens the linear relaxation, possibly improving the solution time.

The rampdown is modeled analogously:
\begin{alignat*}{2}
	& \forallJ,\, k \in \discrange{2}{T}:\quad& \prod(j,k) & \geq \prod(j,k-1) - L \cdot \rampdown(j) \onOff(j,k) - \shutdownRampdown(j) (1 - \onOff(j,k))\\
	&&&+ \min\{\shutdownRampdown(j), \max\{\minProd(j), 0\} + L \cdot \rampdown(j)\} \cdot (1 - \onOff(j,k-1))
\end{alignat*}
Finally, we need to put a limit on the possible maximal production in case of an
imminent shutdown. A shutdown scheduled for period $k+1$ can not be postponed, even in case of a power outage. 
Thus, in case of a shutdown in $k+1$, the unit cannot produce more than $\shutdownRampdown(j)$ in period $k$:
\[\forallJ,\, k \in \discrange{1}{T-1}:\quad \maxPossProd(j,k) \leq \shutdownRampdown(j) (\onOff(j,k) - \onOff(j,k+1)) + \maxProd(j) \onOff(j,k+1)\]
The term $\shutdownRampdown(j) (\onOff(j,k) - \onOff(j,k+1))$ expresses the
desired limit, while the term $\maxProd(j) \onOff(j,k+1)$ keeps the constraint
valid in each of the other three settings of $\onOff(j,k)$ and $\onOff(j,k+1)$.

\subsubsection{Storage (constraints \ref{eq:StorageVariables}, \ref{eq:StorageStart} to \ref{eq:StorageEnd} in subsection \ref{sec:Model})}
\label{sec:Storage}

Storage units are used to even out the demand. The classical examples for storage units are water pumping stations and, to a lesser extent, batteries. We indicate a storage unit by a negative minimal production $\minProd(j)$, meaning that the maximal storage inflow is $-\minProd(j)$. 

\noindent Not being thermal units, the nowadays relevant storage units have very high ramping speeds, which allow us to neglect them, i.e. we assume $\rampup(j) = \rampdown(j) = \shutdownRampdown(j) = \startupRampup(j) = \maxProd(j) + (- \minProd(j))$ for each storage unit $j$.

In our model storage units are characterized by five parameters:
\begin{itemize*}
  \item their storage capacity $\storageCapacity(j)$,
  \item their storage efficiency $\storageEfficiency(j)$,
  \item their constant inflow $\storageInflow(j)$,
  \item and their initial and final storage fill, $\initialStorage(j)$ and $\finalStorage(j)$,
\end{itemize*}
and two new variables: 
\begin{itemize*}
\item the current storage fill $\forallJK: \storage(j, k) \in \set{R}^+$ and
\item the power consumption $\forallJK: \consumption(j, k) \in \set{R}^+$.
\end{itemize*}
 
\noindent The maximal storage inflow and the storage capacity are now easily formulated as
\begin{alignat*}{2}
	& \forallJK: \qquad & \storage(j, k) &\leq \storageCapacity(j)\\
	& \forallJK: \qquad & \consumption(j, k) &\leq \max\{0, - \minProd(j)\}
\end{alignat*}

\noindent The storage fill change from period $k-1$ to $k$ has to account for the prior storage, the stored and consumed energy, and a possible constant inflow (for example a stream feeding a reservoir):
\[\forallJ,\, k \in \discrange{2}{T}: \qquad \storage(j, k) = \storage(j, k-1) + L \cdot (\storageEfficiency(j) \consumption(j, k-1) - \prod(j, k-1) + \storageInflow(j))\]

\noindent The initial and final storage fill parameters $\initialStorage(j)$ and $\finalStorage(j)$ are set outside the model (possibly for raising the total storage when demands and prizes are expected to increase after the considered total time period). They apply to the storage fill at $k = 1$ and $k = T+1$:
\begin{alignat*}{2}
 & \forallJ: \qquad & \storage(j, 1) &= \initialStorage(j)\\
 & \forallJ:        & \finalStorage(j) &= \storage(j, T) + L \cdot (\storageEfficiency(j) \consumption(j, T) - \prod(j, T) + \storageInflow(j))
\end{alignat*}



\subsection{Power Grid Constraints (constraints \ref{eq:DemandAndReserveStart} and \ref{eq:DemandAndReserveEnd} in subsection \ref{sec:Model})}
\label{sec:DemandAndReserve}

All units together have to satisfy the energy demand in each period, given as $\demand(k)$:
\begin{equation*}
	\forallK: \qquad \sum_{j \in J} (\prod(j,k) - \consumption(j, k)) = \demand(k)
\end{equation*}
Moreover, to be failure-tolerant, the units should be able to compensate a power outage, caused for example by a failing unit. This is expressed by the need to keep a given reserve $\reserve(k)$ capacity available:
\begin{equation*}
	\forallK: \qquad \sum_{j \in J} (\maxPossProd(j,k) - \prod(j, k) + \consumption(j, k)) \geq \reserve(k)
\end{equation*}
Typically $\reserve(k)$ is constant over time. However, a dependency e.g. on $\demand(k)$ maybe reasonable, too.


\subsection{Cost Constraints}
\label{sec:CostConstraints}

\subsubsection{Objective Function (constraints \ref{eq:ObjectiveFunction} and \ref{eq:CostVariables} in subsection \ref{sec:Model})}
\label{sec:ObjectiveFunction}

The objective is to minimize the overall costs consisting of production costs
$\prodCost(j,k) \ge 0$ (see subsection \ref{sec:ProductionCosts} below) as well as start-up
and shutdown costs, $\startupCost(j,k), \shutdownCost(j,k) \ge 0$ (see \ref{sec:StartupAndShutdownCosts}):
\[\min \sum_{j \in J} \sum_{k \in K} \prodCost(j,k) + \startupCost(j,k) + \shutdownCost(j,k).\]

\subsubsection{Production Costs (constraint \ref{eq:ProductionCosts} in subsection \ref{sec:Model})}
\label{sec:ProductionCosts}

In \cite{carrion2006computationally} a convex production cost function is assumed and approximated by piecewise affine linear functions. However, since the efficiency of a power unit usually increases with its production level, we assume a concave production cost function. Fortunately, the increase in efficiency is typically quite small, allowing a nearly affine linear approximation
\begin{alignat*}{2}
	& \forallJK: \qquad& \prodCost(j,k) &= A \cdot L \cdot \prod(j,k) + B \cdot L \cdot \onOff(j, k),
\end{alignat*}
with parameters $A$ and $B$. $A$ can be interpreted as the variable cost, incurred for every additional MWh of production, and $B$ as the fixed cost, incurred for running the unit for one hour.

\clearpage\noindent In our model we further split the parameters $A$ and $B$ in two parts, one part which is attributed to the needed fuel, and is thus dependent on the current fuel price, while the other part is attributed to operation and maintenance, and thus fixed. We expect the fuel price to be given as the parameter $\fuelCost(\fuelType(j),k)$, where $\fuelType(j)$\footnote{$\fuelType(j)$ may not contain special characters and spaces, see subsection \ref{sec:InputPeriod}} denotes the fuel type used by unit $j$.

Therefore, we replace $A$ and $B$ by their parts, the time and fuel-type dependend $\varFuel(j)$ and independend $\fixedFuel(j)$ fuel needs and the time and fuel-type dependend $\varCost(j)$ and independend $\fixedCost(j)$ costs. This leads to
\begin{alignat*}{2}
	& \forallJK: \qquad& \prodCost(j,k) &= \left(\varFuel(j) \cdot \fuelCost(\fuelType(j),k) + \varCost(j)\right) L \cdot \prod(j,k)\\
	&                  &                &+ \left(\fixedFuel(j) \cdot \fuelCost(\fuelType(j),k) + \fixedCost(j)\right) L \cdot \onOff(j,k).
\end{alignat*}
Note that due to the affine, non-linear cost function, we do not have a constant production efficiency as modeled in many other papers; instead, the modeled production efficiency increases with the production and is concave, which is typical for thermal units.


\subsubsection{Startup and Shutdown Costs (constraints \ref{eq:StartupAndShutdownCostsStart} and \ref{eq:StartupAndShutdownCostsEnd} in subsection \ref{sec:Model})}
\label{sec:StartupAndShutdownCosts}

Every unit incurs costs when starting up or shutting down, the former increasing with the offline time (e.g. for a thermal unit, the start-up costs are partially attributed to the need for reheating it), the latter constant.

We expect the start-up cost of unit $j$ after an offline time of $t$ periods to be given as $\StartupCost(j,t)$ and the shutdown costs to be given as a single constant $\ShutdownCost(j)$. The start-up and shutdown costs can then be modeled as
\begin{alignat*}{2}
&	 \forallJ,\,k \in \discrange{2}{T}:  \qquad & \shutdownCost(j,k) &\geq \ShutdownCost(j) (\onOff(j,k-1) - \onOff(j,k))\\
&	 \forallJK,\, t \in \discrange{1}{k-1}: \qquad & \startupCost(j,k) &\geq \StartupCost(j,t) \left( \onOff(j,k) - \sum_{n=1}^t \onOff(j,k - n) \right).
\end{alignat*}
Here, the term $\onOff(j,k) - \sum_{n=1}^t \onOff(j,k - n)$ is 1 exactly if unit $j$ starts up in period $k$ and was offline in periods $\discrange{(k - 1)}{(k - t)}$, otherwise it is less or equal to $0$. Thus, the latter constraint is equivalent to
\begin{alignat*}{2}
&	\forallJK,\, t \in \{t' \in \discrange{1}{k - 1}: \forall\ i \in \discrange{1}{t'}: \onOff(j, k - i) = 0\}: \qquad &\startupCost(j,k) &\geq \StartupCost(j,t).
	\intertext{Furthermore, since $\StartupCost(j, t)$ is increasing with $t$, this again is equivalent to}
&	\forallJK,\, \overline{t} = \max\{t' \in \discrange{1}{k - 1}: \forall\ i \in \discrange{1}{t'}: \onOff(j, k - i) = 0\}: \qquad &\startupCost(j,k) &\geq \StartupCost(j,\overline{t}).
\end{alignat*}
Since we minimize the costs, $\startupCost(j,k) = \StartupCost(j,{\overline{t}})$ for start-ups in an optimal solution, as intended.


\subsubsection{Thinning Out the Startup Cost Function}
\label{sec:ThinningOutTheStartupCostFunction}

The number of constraints needed in \cite{carrion2006computationally} to describe the start-up cost function is quite substantial. While all the other constraints amount to about $12 |J||K|$, the model needs about $|J||K|^2$ constraints to model the start-up cost function as discussed in the last sections. Of course it is not necessary to model the start-up cost function for $|K|$ cooldown periods, but only for the first periods, when the cost changes significally. Still, for typical thermal units the relevant timespan is about 2 to 3 days, leading to between $48 |J||K|$ and $72 |J||K|$ constraints if we assume hourly periods, for example.

If we want to reduce the number of constraints further, we have to use a single step for a whole group of cooldown periods, which amounts to assigning the same start-up cost to each of these periods. As an example, let us have a look at a typical start-up cost function with a fixed cost of about $70\%$ and a variable cost of about $30\%$. It is possible to reduce the number of steps from $71$ to $9$, while maintaining a relative error of less than $5\%$:
\begin{center}
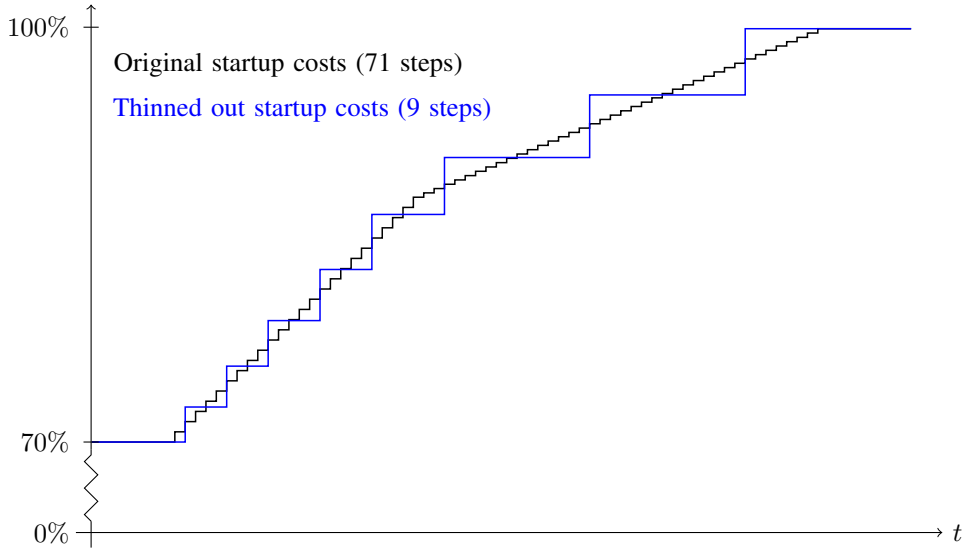

	\captionsetup{type=figure}
	\begin{tikzpicture}[x=1mm,y=1mm]
		\draw[->] (-2,0) -- (112,0);
		\draw (0, -2) -- (0, 1.5);
		\draw decorate [decoration=zigzag] {(0, 1.5) -- (0,11)};
		\draw[->] (0,11) -- (0, 70);
		
		\pgfdeclareimage{function}{\Graphics/StartupCostThinning_func}	
		\pgftext[at=\pgfpoint{0mm}{11.9mm},left,base]{\pgfuseimage{function}};
		
		\node at (114,0) {$t$};
		
		\node[text width=8mm,align=right] at (-7,0) {$0\%$};

		\node[text width=8mm,align=right] at (-7,12) {$70\%$};
		\draw (-1,12) -- (1,12);
		\node[text width=8mm,align=right] at (-7,67) {$100\%$};
		\draw (-1,67) -- (1,67);
		
		\node[text width=120mm,align=left] at (63,62) {Original startup costs (71 steps)};
		\node[text width=120mm,align=left, text=blue] at (63,56) {Thinned out startup costs (9 steps)};
	\end{tikzpicture}

	\captionof{figure}{Approximation of the start-up cost function with a tolerance of $5\%$\\}
	\label{StartupCostThinning}
\end{center}

\noindent Thus, two questions are to be answered: How should we group the periods together to obtain a relative error as small as possible, and which start-up cost should we assign to each of these group?

The second question is answered easily: If we are given a group of periods with cooldown times $t_a$ to $t_b$, it is straight-forward to show that
the minimal relative error\vskip-1ex
\[\text{bestError}(\StartupCost(j,{t_a}), \StartupCost(j,{t_b})) = \frac{\StartupCost(j,{t_b}) - \StartupCost(j,{t_a})}{\StartupCost(j,{t_b}) + \StartupCost(j,{t_a})}\]
\vskip-0.35em\noindent is obtained with the step value\vskip-0.95em
\[\text{bestStep}(\StartupCost(j,{t_a}), \StartupCost(j,{t_b})) = \frac{2\cdot\StartupCost(j,{t_a})\cdot\StartupCost(j,{t_b})}{\StartupCost(j,{t_b}) + \StartupCost(j,{t_a})},\]
and $0$, if $\StartupCost(j,{t_a}) = \StartupCost(j,{t_b}) = 0$.

\clearpage\noindent To answer the first question: let \mosel{STARTUP_TOL} be a parameter giving the maximal error tolerance. Then the periods may be grouped iteratively by starting the first group at $t = 1$, expanding the group as long as its relative error is less than the tolerance, and continuing with the next group on the remaining $t$'s. In pseudo-code, for given \mosel{j}:
\begin{algorithm}
	\DontPrintSemicolon
	\caption{startupCostThinning}
	$t_a, t_b \leftarrow 1$ \tcp*{Start and end index of current group}	
	\While{$t_a \leq |K|$}{
		\tcp{Expand the group as long as the next relative error is less than {\ttfamily STARTUP\_TOL}}
		\While{$t_b + 1 \leq |K|\ \land$ \emph{bestError}$(\StartupCost(j,t_a), \StartupCost(j,t_b+1)) <  $ \emph{\ttfamily STARTUP\_TOL}\ \hbox{}}{
			$t_b \leftarrow t_b + 1$\;
		}
		
		$\StartupCost(j,t_a) \leftarrow \text{bestStep}(\StartupCost(j,t_a), \StartupCost(j,t_b))$
		\tcp*{Calculate optimal step value}
		$\forall\ t \in \discrange{(t_a\!+\!1)}{t_b}: \StartupCost(j,t) \leftarrow 0$
		\tcp*{Delete other step values}
		$t_a, t_b \leftarrow t_b + 1$ \tcp*{Continue with next group}
	}
\end{algorithm}                

\noindent It can be proven that the number of groups produced by this algorithm is the minimal number of groups needed to fulfill the tolerance requirements.


\subsection{Dealing with Infeasibilities}
\label{sec:DealingWithInfeasibilities}

The main cause of infeasibilities is an excessive demand or reserve capacity,
\[\exists\ k: \sum_{j \in J} \maxProd(j, k) < \demand(j) + \reserve(j).\]
In most cases the demand not satisfiable by the units within the model is covered by energy sources which are not part of the model (e.g. by power imports from external markets or backup units). Therefore it is advisable to soften the demand and reserve constraints by measuring the underproduction and underreserve and applying a penalty to it. The softened demand and reserve constraints thus are:
\begin{alignat*}{2}
	& \forallK: \qquad & \sum_{j \in J} (\prod(j,k) - \consumption(j,k)) + \underProduction(k) &= \demand(k)\\
	& \forallK:        & \sum_{j \in J} (\maxPossProd(j,k) - \prod(j,k) + \consumption(j,k)) + \underReserve(k) &\geq \reserve(k)
\end{alignat*}

\noindent Given the two parameters $\underProductionPenalty$ (penalty factor for underproduction) and $\underReservePenalty$ (penalty factor for underreserve), we replace the objective function by
\[\min \sum_{j \in J} \sum_{k \in K} \prodCost(j,k) + \startupCost(j,k) + \shutdownCost(j,k) + \sum_{k \in K} \underProductionPenalty \cdot \underProduction(k) + \underReservePenalty \cdot \underReserve(k).\]

\clearpage\noindent
The second most common cause is a demand profile too volatile to be satisfied by the units, or in other words, units too slow to follow the demand profile. This infeasibility also originates from energy sources not present in the model. Although the softening against underproduction is enough to prevent infeasibilities, one may want to penalize overproduction less (in real world, overproduction may for example be dumped on external markets). So we obtain
\begin{alignat*}{2}
	& \forallK: \qquad & \sum_{j \in J} (\prod(j,k) - \consumption(j,k)) + \underProduction(k) - \overProduction(k) &= \demand(k)\\
	& \forallK:        & \sum_{j \in J} (\maxPossProd(j,k) - \prod(j,k) + \consumption(j,k)) + \underReserve(k) &\geq \reserve(k)
\end{alignat*}
\noindent and, given the parameter $\overProductionPenalty$ (penalty factor for overproduction),
\[\min \sum_{j \in J} \sum_{k \in K} \prodCost(j,k) + \startupCost(j,k) + \shutdownCost(j,k) + \sum_{k \in K} \underProductionPenalty \cdot \underProduction(k) + \underReservePenalty \cdot \underReserve(k) + \overProductionPenalty \cdot \overProduction(k).\]

\noindent There are a few other inconsistencies in the input data which could lead to infeasibilities or unwanted behavior, namely
\begin{itemize}
  \item simultaneous positive initial down- and uptime:\\
  	$\exists\ j \in J: \initialOnline(j) > 0\ $ and $\ \initialOffline(j) > 0$,
  \item impossible production and ramping limits:\\
  	$\exists\ j \in J: \minProd(j) > \maxProd(j)\ $ or $\ \minProd(j) > \startupRampup(j)\ $ or $\ \minProd(j) > \shutdownRampdown(j)$,
  \item decreasing start-up costs:\\
  	$\exists\ j \in J,\, k \in K: \StartupCost(j, k) > \StartupCost(j, k + 1)$, and
  \item non fulfillable storage constraints:\\
    $\exists\ j \in J: \storageInflow(j) > \maxProd(j)$\\
    $\exists\ j \in J: \storageEfficiency(j) \notin [0, 1]$\\
    $\exists\ j \in J: \initialStorage(j) > \storageCapacity(j)\ $ or $\ \finalStorage(j) > \storageCapacity(j)$\\
    $\exists\ j \in J: \initialStorage(j) + L \cdot T \cdot (- \storageEfficiency(j) \minProd(j) + \storageInflow(j)) < \finalStorage(j)\ $ or $\ \initialStorage(j) + L \cdot T \cdot (- \maxProd(j) + \storageInflow(j)) > \finalStorage(j)$
\end{itemize}
All those are errors in the input data, which have to be avoided for meaningful optimization.

\clearpage
\subsection{Parameters summary}
\label{sec:Parameters}

\subsubsection{General parameters}
\begin{table}[H]
\begin{tabular}{l|c|c|l}
Name & Domain & Scale unit & Description\\
\hline
$\periodLength$ & $\set{R}^+$ & h & Period length\\
$\underProductionPenalty$ & $\set{R}^+$ & cost/MWh & Penalty for underproduction\\
$\underReservePenalty$ & $\set{R}^+$ & cost/MWh & Penalty for underreserve\\
$\overProductionPenalty$ & $\set{R}^+$ & cost/MWh & Penalty for overproduction\\
STARTUP\_TOL & [0,1] & & Startup cost error tolerance
\end{tabular}
\end{table}

\subsubsection{Parameters for each period}

\begin{table}[H]
\begin{tabular}{l|c|c|l}
Name & Domain & Scale unit & Description\\
\hline
$\demand(k)$ & $\set{R}^+$ & MW & Power demand\\
$\reserve(k)$ & $\set{R}^+$ & MW & Needed reserve capacity\\
$\fuelCost(f,k)$ & $\set{R}^+$ & cost/MWh & Cost of a MWhs worth of a fuel type\\
\end{tabular}
\end{table}

\subsubsection{Parameters for each power unit}

\begin{table}[H]
\begin{tabular}{l|c|c|l}
Name & Domain & Scale unit & Description\\
\hline
$\minUptime(j)$ & $\discrange{1}{T}$ & periods & Minimal uptime after starting up\\
$\minDowntime(j)$ & $\discrange{1}{T}$ & periods & Minimal downtime after shutting down\\
$\initialOnline(j)$ & $\discrange{0}{T}$ & periods & Initial uptime (unit is online in periods $1, \dots, \initialOnline(j)$) \\
$\initialOffline(j)$ & $\discrange{0}{T}$ & periods & Initial downtime (unit is offline in periods $1, \dots, \initialOffline(j)$) \\
$\minProd(j)$ & $\set{R}$ & MW & Minimal production in online state\\
$\maxProd(j)$ & $\set{R}^+$ & MW & Maximal production in online state\\
$\rampup(j)$ & $\set{R}^+$ & MW/h & Maximal production increase per hour, in online state\\
$\rampdown(j)$ & $\set{R}^+$ & MW/h & Maximal production decrease per hour, in online state\\
$\startupRampup(j)$ & $\set{R}^+$ & MW & Maximal production increase at start-up\\
$\shutdownRampdown(j)$ & $\set{R}^+$ & MW & Maximal production decrease on shutdown\\
$\storageCapacity(j)$ & $\set{R}^+$ & MWh & Storage capacity (Zero for most units)\\
$\storageEfficiency(j)$ & $[0, 1]$ & MWh/MWh & Storage efficiency (Zero for most units)\\
$\storageInflow(j)$ & $\set{R}^+$ & MW & Storage inflow (Zero for most units)\\
$\initialStorage(j)$ & $\set{R}^+$ & MWh & Initial storage fill\\
$\finalStorage(j)$ & $\set{R}^+$ & MWh & Final storage fill\\
\hline
$\fuelType(j)$ & & fuel & Used fuel type \footnotemark[1]\\
$\varFuel(j)$ & $\set{R}^+$ & MWh/MWh & Fuel need increase per MW of production\\
$\fixedFuel(j)$ & $\set{R}^+$ & MWh/h & Fixed fuel need in online state\\
$\varCost(j)$ & $\set{R}^+$ & cost/MWh & Cost increase per MWh of production\\
$\fixedCost(j)$ & $\set{R}^+$ & cost/h & Fixed cost in online state\\
$\StartupCost(j,t)$ & $\set{R}^+$ & cost & Startup cost after $t$ offline periods\\
$\ShutdownCost(j)$ & $\set{R}^+$ & cost & Shutdown cost\\
\end{tabular}
\label{tbl:UnitParameters}
\end{table}
\footnotetext[1]{$\fuelType(j)$ may not contain special characters and spaces, see subsection \ref{sec:InputPeriod}}

\clearpage
\subsection{Model summary}
\label{sec:Model}

%
%
\noindent Objective function (discussion in \ref{sec:ObjectiveFunction} and \ref{sec:DealingWithInfeasibilities}, implementation in \ref{sec:ImplObjectiveFunction}):
\begin{equation}
\min \sum_{j \in J} \sum_{k \in K} \prodCost(j,k) + \startupCost(j,k) + \shutdownCost(j,k) + \sum_{k \in K} \underProductionPenalty \cdot \underProduction(k) + \underReservePenalty \cdot \underReserve(k) + \overProductionPenalty \cdot \overProduction(k) \label{eq:ObjectiveFunction}
\end{equation}
%
%
Physical unit variables (discussion in \ref{sec:UnitVariables} and \ref{sec:Storage}, implementation in \ref{sec:ImplVariables}):
\begin{alignat}{2}
	& \forallJK: \qquad & \onOff(j,k) &\in \{0,1\} \label{eq:UnitVariablesStart}\\
	& \forallJK:        & \prod(j,k), \maxPossProd(j,k)  &\geq 0 \label{eq:UnitVariablesEnd}\\
    & \forallJK:        & \storage(j, k), \consumption(j, k) &\geq 0 \label{eq:StorageVariables}
\intertext{Cost variables (discussion in \ref{sec:CostConstraints}, implementation in \ref{sec:ImplVariables})}
    & \forallJK:        & \prodCost(j,k), \shutdownCost(j,k), \startupCost(j,k) &\geq 0 \label{eq:CostVariables}
\intertext{Period variables (discussion in \ref{sec:DemandAndReserve} and \ref{sec:DealingWithInfeasibilities}, implementation in \ref{sec:ImplVariables}):}
    & \forallK:         & \underProduction(k), \overProduction(k), \underReserve(k) &\geq 0
\end{alignat}
%
%
Minimal up- and downtime (discussion in \ref{sec:MinimalUpAndDowntime}, implementation in \ref{sec:ImplPhysicalConstraintsUnits}):
\begin{alignat}{2}
	& \forallJ,\, k \in \discrange{1}{\initialOnline(j)}: \qquad & \onOff(j,k) &= 1 \label{eq:MinimalUpAndDowntimeStart}\\
	& \forallJ,\, k \in \discrange{1}{\initialOffline(j)}:        & \onOff(j,k) &= 0
\end{alignat}\vskip-3\intextsep\begin{alignat}{2}
\forallJ,\, k \in \discrange{(\initialOnline(j) + 2)}{T},\, i \in \discrange{1}{\minUptime(j) - 1} \cap \discrange{1}{T - k}:&\qquad\onOff(j,k+i) \geq& \onOff(j,k) - \onOff(j,k-1)\phantom{)}\\
\forallJ,\, k \in \discrange{(\initialOffline(j) + 2)}{T},\, i \in \discrange{1}{\minDowntime(j) - 1} \cap \discrange{1}{T - k}:&\qquad \onOff(j,k+i) \leq&\ 1 - (\onOff(j,k) - \onOff(j,k-1))
\label{eq:MinimalUpAndDowntimeEnd}
\end{alignat}
Minimal and maximal production (discussion in \ref{sec:MinimalAndMaximalProduction}, implementation in \ref{sec:ImplPhysicalConstraintsUnits}):
\begin{equation}
\forallJK: \qquad \minProd(j) \onOff(j,k) \leq \prod(j,k) \leq \maxPossProd(j,k) \leq \maxProd(j) \onOff(j,k) \label{eq:MinimalAndMaximalProduction}
\end{equation}
Ramping constraints (discussion in \ref{sec:Ramping}, implementation in \ref{sec:ImplPhysicalConstraintsUnits}):
\begin{alignat}{2}
	& \forallJ,\, k \in \discrange{2}{T}: \qquad & \maxPossProd(j,k) & \leq \prod(j,k-1) + L \cdot \rampup(j) \onOff(j,k-1) + \startupRampup(j) (1 - \onOff(j,k-1)) \label{eq:RampingStart}\\
	&&&\quad - \min\{\startupRampup(j), \minProd(j) + L \cdot \rampup(j) \} \cdot (1 - \onOff(j,k)) \nonumber\\
	& \forallJ,\, k \in \discrange{2}{T}: & \prod(j,k) & \geq \prod(j,k-1) - L \cdot \rampdown(j) \onOff(j,k) - \shutdownRampdown(j) (1 - \onOff(j,k))\\
	&&&\quad + \min\{\shutdownRampdown(j), \minProd(j) + L \cdot \rampdown(j)\} \cdot (1 - \onOff(j,k-1)) \nonumber\\
	& \forallJ,\, k \in \discrange{1}{T-1}: & \maxPossProd(j,k) &\leq \maxProd(j) \onOff(j,k+1) + \shutdownRampdown(j) (\onOff(j,k) - \onOff(j,k+1)) \label{eq:RampingEnd}
\end{alignat}
Storage constraints (discussion in \ref{sec:Storage}, implementation in \ref{sec:ImplPhysicalConstraintsUnits}):
\begin{alignat}{2}
 & \forallJK: \qquad & \storage(j, k) &\leq \storageCapacity(j)\label{eq:StorageStart}\\
 & \forallJK:        & \consumption(j, k) & \leq \max\{0, - \minProd(j)\}\\
 & \forallJ,\, k \in \discrange{2}{T}: \qquad & \storage(j, k) &= \storage(j, k-1) + L \cdot (\storageEfficiency(j) \consumption(j, k-1) - \prod(j, k-1) + \storageInflow(j))\\
 & \forallJ:         & \storage(j, 1) &= \initialStorage(j)\\
 & \forallJ:         & \finalStorage(j) &= \storage(j, T) + L \cdot (\storageEfficiency(j) \consumption(j, T) - \prod(j, T) + \storageInflow(j)) \label{eq:StorageEnd}
\end{alignat}
Demand and reserve (discussion in \ref{sec:DemandAndReserve} and \ref{sec:DealingWithInfeasibilities}, implementation in \ref{sec:ImplDemandAndReserve}):
\begin{alignat}{2}
	& \forallK: \qquad & \sum_{j \in J} (\prod(j,k) - \consumption(j,k)) + \underProduction(k) - \overProduction(k) &= \demand(k) \label{eq:DemandAndReserveStart}\\
	& \forallK:        & \sum_{j \in J} (\maxPossProd(j,k) - \prod(j,k) + \consumption(j,k)) + \underReserve(k) & \geq \reserve(k) \label{eq:DemandAndReserveEnd}
\end{alignat}
Production costs (discussion in \ref{sec:ProductionCosts}, implementation in \ref{sec:ImplProductionCosts}):
\begin{alignat}{2}
	& \forallJK: \qquad& \prodCost(j,k) &= \left(\varFuel(j) \cdot \fuelCost(\fuelType(j),k) + \varCost(j) \right) L \cdot \prod(j,k) \label{eq:ProductionCosts}\\
	&                                  &&\quad + \left(\fixedFuel(j) \cdot \fuelCost(\fuelType(j),k) + \fixedCost(j) \right) L \cdot \onOff(j,k)\nonumber
\intertext{Startup and shutdown costs (discussion in \ref{sec:StartupAndShutdownCosts} and \ref{sec:ThinningOutTheStartupCostFunction}, implementation in \ref{sec:ImplStartupAndShutdownCosts}):}
	&	 \forallJ,\,k \in \discrange{2}{T}:  \qquad & \shutdownCost(j,k) &\geq \ShutdownCost(j) (\onOff(j,k-1) - \onOff(j,k))\label{eq:StartupAndShutdownCostsStart}\\
	&	 \forallJK,\, t \in \discrange{1}{k-1}: \qquad & \startupCost(j,k) &\geq \StartupCost(j,t) \left( \onOff(j,k) - \sum_{n=1}^t \onOff(j,k - n) \right)\label{eq:StartupAndShutdownCostsEnd}
\end{alignat}

\clearpage

\section{Implementing the Model}
\label{sec:Implementation}

In this section, we show how to implement the model discussed in the last section in Mosel \cite{mosel}, i.e. how the model composed by parameters, variables and constraints is stated and solved. The how-to on the surrounding tasks, data input and output, is given in the following sections.

For the whole implementation, we assume the following \mosel{options} to be used:
\begin{Mosel}{}
options noimplicit, explterm, keepassert;
\end{Mosel}
The option \mosel{noimplicit} disallows implicit variable creation, and thus forces the programmer to declare each variable. Mosel's automatic line termination requires us to wrap multi-line formulas after an operator, therefore we deactivate it with \mosel{explterm}. We check the input data with \mosel{assert()}, which by default is only active in debug mode; \mosel{keepassert} enables it in all modes.

\subsection{Parameters (subsection \ref{sec:Parameters})}
\label{sec:ImplParameters}

\noindent The fitting datatype is clear from the parameters table: real parameters are represented as \mosel{real}, integer parameters as \mosel{integer} and fuel types as \mosel{string}. The index sets of the parameters are either \mosel{J}, \mosel{K}, the fuel types or a combination of two of these sets. For performance reasons, the order of the index sets should be the same as in the \mosel{forall} loops defined over these sets. Here, this means that the set \mosel{K} always comes last.

Based on these design decisions, we declare the parameters as
\begin{Mosel}{}
declarations
	J: set of integer;				  ! Set of units
	K = 1..T;						  ! Set of periods
	Fuels: set of string;			  ! Available fuel types (index set for CF)
	UPP, URP, OPP: real;			  ! Penalties [cost/MWh]

	D, R: array(K) of real;			  ! Demand and reserve [MW]
	FC: dynamic array(Fuels) of
			array(K) of real;         ! Cost of a fuel type [cost/MWh]

	UT, DT: array(J) of integer;	  ! Minimal up-/downtime [periods]
	IUT, IDT: array(J) of integer;	  ! Initial up-/downtime [periods]
	P_min, P_max: array(J) of integer;! Minimal and maximal production [MW]
	RU, RD: array(J) of integer;	  ! Maximal upwards/downwards ramping [MW/h]
	SU, SD: array(J) of integer;	  ! Startup/shutdown ramping [MW]
	SC: array(J) of real;			  ! Storage capacity [MWh]
	SE: array(J) of real;			  ! Storage efficiency [MWh/MWh]
	SIF: array(J) of real;			  ! Storage inflow [MW]
	SI, SF: array(J) of real;		  ! Initial and final storage fill [MWh]
	F: array(J) of string;			  ! Fuel type of a unit [Fuels]
	FA, PA: array(J) of real;		  ! Variable prod. costs [MW/MW],[cost/MW]
	FB, PB: array(J) of real;		  ! Fixed production costs [MW],[cost]
	CU: dynamic array(J, range) of real;! Startup costs [cost]
	CD: array(J) of real;			  ! Shutdown costs [cost]
end-declarations
\end{Mosel}

\noindent The parameters \mosel{T} (number of periods) and \mosel{L} (period length) are declared in the data input section (\ref{sec:DataInput}), since they need to be known for selecting the data to be read.

The size of \mosel{CU}'s second index set \mosel{range} is unknown, since we do not know in advance how many different start-up costs are available for each unit. We therefore declare this array as \mosel{dynamic}.

\mosel{FC} has to be a dynamic array too, since the set of fuels is not known in advance. Usually, we would declare such an array as
\begin{Mosel}{}
	FC: dynamic array(Fuels, K) of real;
\end{Mosel}
However, using the more specific declaration
\begin{Mosel}{}
	FC: dynamic array(Fuels) of array(K) of real;
\end{Mosel}
has an important advantage: This way we may pass individual columns of the matrix, the subarrays \mosel{FC(f)} to the SQL functions, see subsection \ref{sec:InputPeriod}. Conveniently, the elements of \mosel{FC} can still be accessed as \mosel{FC(f, k)}.


\subsection{Variables}
\label{sec:ImplVariables}

\noindent Again, the index sets of the variables are determined from the used indices. In Mosel, the usual datatype for variables is \mosel{mpvar} (mathematical programming decision variable). The only exception is the objective function, which as a linear function is declared as \mosel{linctr} (linear constraint, also used for linear functions).

As the number of storage units is usually quite small, we declare the storage and consumption variable arrays as \mosel{dynamic} and create the individual variables only for them, thus saving variables.
\begin{Mosel}{}
declarations
	p: array(J, K) of mpvar;			! Production [MW]
	p_max: array(J, K) of mpvar;		! Maximal possible production [MW]
	v: array(J, K) of mpvar;			! On-Off state [0/1]
	s: dynamic array(J, K) of mpvar;	! Storage [MWh]
	c: dynamic array(J, K) of mpvar;	! Consumption [MW]
	cp: array(J, K) of mpvar;			! Production costs [cost]
	cu, cd: array(J, K) of mpvar;		! Startup/shutdown costs [cost]
	overallCosts: linctr;				! Overall costs [cost]
	p_under, p_over: array(K) of mpvar;	! Under-/overproduction [MW]
	r_under: array(K) of mpvar;			! Underreserve [MW]
end-declarations

! Create storage and consumption variables for storage units
forall(j in J, k in K | P_min(j) < 0) do
	create(s(j, k));
	create(c(j, k));
end-do
\end{Mosel}

\noindent By default, these variables are constrained to take positive real values. This fits all variables except the on-off-state \mosel{v}, which should be binary:
\begin{Mosel}{}
forall(j in J, k in K) v(j, k) is_binary;
\end{Mosel}


\subsection{Physical Constraints of the Units (constraints \ref{eq:MinimalUpAndDowntimeStart} to \ref{eq:RampingEnd})}
\label{sec:ImplPhysicalConstraintsUnits}

The Mosel equivalent of the {\bf minimal up- and downtime} constraints (\ref{eq:MinimalUpAndDowntimeStart} to \ref{eq:MinimalUpAndDowntimeEnd}) is canonical:
\begin{Mosel}{}
forall(j in J, k in 1..IUT(j)) v(j, k) = 1;
forall(j in J, k in 1..IDT(j)) v(j, k) = 0;
\end{Mosel}
for the initial up- and downtime, and
\begin{Mosel}{}
forall(j in J, k in IUT(j)+2 .. T, i in 1..minlist(UT(j)-1, T-k))
	v(j, k + i) >= v(j, k) - v(j, k-1);
forall(j in J, k in IDT(j)+2 .. T, i in 1..minlist(DT(j)-1, T-k))
	v(j, k + i) <= 1 - (v(j, k-1) - v(j, k));
\end{Mosel}
for the interperiod minimal up- and downtime.

\noindent The jointed constraints for {\bf minimal and maximal production} (\ref{eq:MinimalAndMaximalProduction}) have to be separated:
\begin{Mosel}{}
forall(j in J, k in K) do
	P_min(j) * v(j, k) <= p(j, k);
	p(j, k) <= p_max(j, k);
    p_max(j, k) <= P_max(j) * v(j, k);
end-do
\end{Mosel}

\noindent The {\bf ramping limits} (\ref{eq:RampingStart} to \ref{eq:RampingEnd}) are straight-forward to implement:
\begin{Mosel}{}
forall(j in J, k in 2..T) do
	p_max(j, k) <= p(j, k-1) + L * RU(j) * v(j, k-1) + SU(j) * (1 - v(j, k-1))
	               - minlist(SU(j), P_min(j) + L * RU(j)) * (1 - v(j, k));
	p(j, k)     >= p(j, k-1) - L * RD(j) * v(j, k) - SD(j) * (1 - v(j, k))
	               + minlist(SD(j), P_min(j) + L * RD(j)) * (1 - v(j, k-1));
end-do
forall(j in J, k in 1..T-1) do
	p_max(j, k) <= P_max(j) * v(j, k+1) + SD(j) * (v(j, k) - v(j, k+1));
end-do
\end{Mosel}

\noindent Same as with the storage variables, the {\bf storage constraints}  (\ref{eq:StorageStart} to \ref{eq:StorageEnd}) are implemented only for storage units.
\begin{Mosel}{}
forall(j in J | P_min(j) < 0) do
	forall(k in K) do
		s(j, k) <= SC(j);
		c(j, k) <= maxlist(0, -P_min(j));
	end-do
	forall(k in 2..T) do
		s(j, k) = s(j, k-1) + L * (SE(j)*c(j, k-1) - p(j, k-1) + SIF(j));
	end-do
	SI(j) = s(j, 1);
	SF(j) = s(j, T) + L * (SE(j)*c(j, T) - p(j, T) + SIF(j));
end-do
\end{Mosel}


\clearpage
\subsection{Power Grid Constraints (constraints \ref{eq:DemandAndReserveStart} and \ref{eq:DemandAndReserveEnd})}
\label{sec:ImplDemandAndReserve}

Since the consumption variables are only created for storage units, we use \mosel{exists} to sum only the existing variables.
\begin{Mosel}{}
forall(k in K) do
    sum(j in J) p(j, k) - sum(j in J | exists(c(j, k))) c(j, k)
        + p_under(k) - p_over(k) = D(k);
    sum(j in J) (p_max(j, k) - p(j, k)) - sum(j in J | exists(c(j, k))) c(j, k)
        + r_under(k) >= R(k);
end-do
\end{Mosel}


\subsection{Production Cost Constraints (constraint \ref{eq:ProductionCosts})}
\label{sec:ImplProductionCosts}

The model could be reduced by replacing every use of $\prodCost(j, k)$ by the right-hand side of constraints \ref{eq:ProductionCosts}, and thus removing the variables $\prodCost(j, k)$  completely. Fortunately, the Xpress Optimizer automatically applies this reduction at the presolve stage and enables us to use the canonical implementation for better readability:
\begin{Mosel}{}
forall(j in J, k in K) cp(j, k) = (FA(j) * FC(F(j), k) + PA(j)) * p(j, k)
                                + (FB(j) * FC(F(j), k) + PB(j)) * v(j, k);
\end{Mosel}


\subsection{Startup and Shutdown Costs (constraints \ref{eq:StartupAndShutdownCostsStart} and \ref{eq:StartupAndShutdownCostsEnd})}
\label{sec:ImplStartupAndShutdownCosts}

The start-up and shutdown constraints can be implemented in Mosel as
\begin{Mosel}{}
forall(j in J, k in 2..T) cd(j, k) >= CD(j) * (v(j, k-1) - v(j, k));
forall(j in J, k in K, t in 1..k-1 | exists(CU(j, t)))
	cu(j, k) >= CU(j, t) * (v(j, k) - sum(n in 1..t) v(j, k - n));
\end{Mosel}
The \mosel{exists} operator is again used to enumerate only relevant indices. Especially the number of constraints depending on $t$ may be further reduced from considering only a subset (compare subsection \ref{sec:ThinningOutTheStartupCostFunction}), accepting a loss in accuracy of the start-up costs.

\subsubsection{Thinning Out the Startup Cost Function}

The pseudo-code of the start-up cost thinning can be translated one-to-one into Mosel code:                                       
\begin{Mosel}{}
function bestError(CU_j_ta: real, CU_j_tb: real): real
	if(CU_j_ta = 0 and CU_j_tb = 0) then
		returned := 0;
	else
		returned := abs(CU_j_ta - CU_j_tb)/(CU_j_ta + CU_j_tb);
	end-if
end-function
\end{Mosel}
\clearpage
\begin{Mosel}{}
function bestStep(CU_j_ta: real, CU_j_tb: real): real
	if(CU_j_ta = 0 and CU_j_tb = 0) then
		returned := 0;
	else
		returned := 2 * CU_j_ta * CU_j_tb / (CU_j_ta + CU_j_tb);
	end-if
end-function
\end{Mosel}
\begin{Mosel}{}
procedure startupCostThinning(j: integer)
	declarations
		t_a, t_b: integer;
	end-declarations
	! Start the first group
	t_a := 1;
	t_b := 1;
	repeat
		! Expand the group as long as the next
		! relative error is less than STARTUP_TOL
		while(exists(CU(j,t_b + 1)) and
		      bestError(CU(j,t_a),CU(j,t_b+1)) < STARTUP_TOL) do
			t_b := t_b+1;
		end-do
		! Calculate optimal step value
		CU(j, t_a) := bestStep(CU(j, t_a), CU(j, t_b));
		! Delete other step values
		forall(t in t_a+1..t_b) do
			delcell(CU(j, t));
		end-do
		! Continue with next group
		t_a := t_b + 1;
		t_b := t_b + 1;
	until(not exists(CU(j, t_a)));
end-procedure
\end{Mosel}
Of course, the thinning has to be applied to every unit $j \in J$:
\begin{Mosel}{}
forall(j in J) do
	startupCostThinning(j);
end-do
\end{Mosel}


\subsection{Objective Function (constraint \ref{eq:ObjectiveFunction})}
\label{sec:ImplObjectiveFunction}

First, we have to define the cost function
\begin{Mosel}{}
overallCosts := sum(j in J, k in K) (cp(j, k) + cd(j, k) + cu(j, k))
              + sum(k in K) (UPP*p_under(k) + URP*r_under(k) + OPP*p_over(k));
\end{Mosel}
and then call \mosel{minimize}:\\
\mosel{minimize(overallCosts);}

\subsection{Dealing with Infeasibilities (subsection \ref{sec:DealingWithInfeasibilities})}
\label{sec:ImplDealingWithInfeasbilities}

The penalties have already been incorporated in the power grid constraints and the objective function. To detect errors in the input data, we use Mosel's \mosel{assert} procedure, which stops the program if a condition is not met.
\begin{Mosel}{}
assert(and(j in J) IUT(j) in 0..T, "Initial uptime out of range!");
assert(and(j in J) IDT(j) in 0..T, "Initial uptime out of range!");
assert(and(j in J) IUT(j)*IDT(j) = 0, "Simultaneous initial down- and uptime!");
assert(and(j in J) UT(j) in 1..T, "Initial uptime out of range!");
assert(and(j in J) DT(j) in 1..T, "Initial downtime out of range!");
assert(and(j in J) P_min(j) <= P_max(j), "Impossible production limits!");
assert(and(j in J) P_min(j) <= SU(j), "Some unit is not able to start up!");
assert(and(j in J) P_min(j) <= SD(j), "Some unit is not able to shutdown!");
declarations
	lastCU: real;
end-declarations
forall(j in J) do
	lastCU := 0;
	forall(k in K | exists(CU(j, k))) do
		assert(lastCU <= CU(j, k),
		       "The start-up costs are not monotonically increasing!");
		lastCU := CU(j, k);
	end-do
end-do
assert(and(j in J) SIF(j) <= P_max(j), "Storage inflow leads to overcapacity!");
assert(and(j in J) (0 <= SE(j) and SE(j) <= 1), "Invalid storage efficiency!");
assert(and(j in J) SI(j) <= SC(j), "Invalid initial storage fill!");
assert(and(j in J) SF(j) <= SC(j), "Invalid final storage fill!");
assert(and(j in J) SI(j) + L*T*(SE(j)*maxlist(0, -P_min(j)) + SIF(j)) >= SF(j),
       "Some storage constraints are not fulfillable!");
assert(and(j in J	) SI(j) + L*T*(-P_max(j) + SIF(j)) <= SF(j),
       "Some storage constraints are not fulfillable!");
\end{Mosel}

\noindent As noted at the beginning of this section (see \ref{sec:Implementation}), we need to use the option \mosel{keepassert} to activate \mosel{assert} outside of debug mode.
\clearpage

\section{Data Input}
\label{sec:DataInput}


\subsection{General Parameters}

We use the two parameters \mosel{UNITS_SOURCE} and \mosel{PERIODS_SOURCE} to store the sources for the unit parameters and the period data. Depending on whether one wants to access a database or a spreadsheet, these parameters contain the name of a database table or a spreadsheet.


Once we know where to get the period data, we have to know a start time \mosel{START} and the number of periods \mosel{T}. Since the datatype \mosel{datetime} can not be used for parameters, \mosel{START} will be a \mosel{string}, formatted with the standard SQL format YYYY-MM-DD HH:MM:SS (example: June 30th, 2007 at 7:12 PM is written as 2007-06-30 19:12:00). The parameter \mosel{PERIOD_LENGTH} should be given in seconds.

Summarising, the used parameters are:
\begin{Mosel}{}
parameters
	UNITS_SOURCE = "Units.xls";		! Source of the unit parameters
	PERIODS_SOURCE = "Periods.xls";	! Source of the period data
	DSN = "Server=?;Database=?;UID=?;PWD=?;"; ! Example DSN
	START = "2009-05-11 00:00:00";	! Start of the modeled timespan
	T = 168;						! Number of modeled periods
	L = 1;							! Period length [h]
	STARTUP_TOL = 0.05;				! Tolerance for modeling of start-up cost
end-parameters
\end{Mosel}
Important: parameters have to be declared with a default value to define their data type.


\subsection{Unit Parameters}

The parameters describing the power units are too many to be read from the command line, so we have
to read them from a different data source. Mosel supports database access, including Excel spreadsheets, over ODBC.

Now, databases are superior to spreadsheets in performance, reliability and
scalability, making them a good choice for a Mosel model used in day-to-day operations.
Excel spreadsheets on the other hand are easy to create and modify, making
them a good choice for the development and experimental stage.

Fortunately, the Mosel \mosel{initializations} block hides most
of the differences between real databases and Excel spreadsheets, allowing us to
switch between them with little effort. For details on how to use the \mosel{initializations} block,
please refer to the FICO Whitepaper \quote{Using ODBC and other database
interfaces with Mosel}, which should reside in the Xpress installation directory at
{\ttfamily .../XpressMP/docs/\linebreak[0]mosel/mosel\_odbc/moselodbc.pdf} \cite{heipcke2009moselodbc}.

We expect two tables in our database,
\begin{itemize}
  \item one table with the unit parameters with index set \mosel{J} (all except start-up costs),
  \item one table with the start-up costs with index sets \mosel{J} and \mosel{K}.
\end{itemize}
The name of the first table should be given as the parameter \mosel{UNITS_SOURCE}, whereas the second table should have the same name with an appended \quote{\_CU} (in reference to the start-up costs name $\StartupCost(j,t)$). The names of the index and of the columns should be the same as in our model.

The Excel spreadsheet is set up in the same way, except that we store each set of units in their own spreadsheet file. Thus the \mosel{UNITS_SOURCE} parameter now denotes the spreadsheet file, while the two sheets are always called \quote{Units} and \quote{Units\_CU}.

Finally, we need the ODBC data source name (DSN), also called \quote{Connection String}, which differs with the brand of database server used. We expect the DSN to be given as \mosel{DSN}.
Connection strings can be found in the database's manual. A collection of connection strings for popular databases is also available at \url{www.connectionstrings.com}.

Now, we can setup the parameters needed by \mosel{initializations}:
\begin{Mosel}{}
declarations
	unitsDSN: string;		! ODBC data source name
	unitsTable: string;		! Table containing the unit parameter
	unitsCUTable: string;	! Table containing the start-up costs
end-declarations
if(isExcelDocument(UNITS_SOURCE)) then
	unitsDSN := "mmodbc.excel:skiph;" + UNITS_SOURCE;
	unitsTable := "Units";
else
	unitsDSN := "mmodbc.odbc:" + DSN;
	unitsTable := UNITS_SOURCE;
end-if
unitsCUTable := unitsTable + "_CU";
\end{Mosel}
The function \mosel{isExcelDocument} used here is listed in Appendix \ref{sec:ExcelFunctions}, and works by checking the file extension.

In the \mosel{initializations} block, we read all unit parameters with index set $J$ from the first table\vskip0pt
\begin{Mosel}{}
initializations from unitsDSN
   [UT,DT,IUT,IDT,P_min,P_max,RU,RD,SU,SD,SC,SE,SIF,SI,SF,F,FA,FB,PA,PB,CD]
   as unitsTable +
   "(j,UT,DT,IUT,IDT,P_min,P_max,RU,RD,SU,SD,SC,SE,SIF,SI,SF,F,FA,FB,PA,PB,CD)";
\end{Mosel}
and the start-up costs from the second table:
\begin{Mosel}{}
	CU as unitsCUTable + "(j,k,CU)";
end-initializations
\end{Mosel}
Mosel automatically uses the columns \mosel{j} and \mosel{k} as the unit and period indices. Once the units have been read, we are able to define the \mosel{Fuels} set:
\begin{Mosel}{}
Fuels := union(j in J) {F(j)};
finalize(Fuels);
\end{Mosel}


\subsection{Period Data}
\label{sec:InputPeriod}

Since we only want to get data for the periods in our timespan, not the full dataset, and we (possibly) do not know all used fuel types, the number of arrays to be read must be variable, which is not possible within one \mosel{initializations}. An elegant solution to this problem is the direct access to the database using Mosel's SQL functions.

The period data is stored in a single table which holds the demand \mosel{D}, the reserve \mosel{R} and the fuel costs. For every used fuel, we expect the apposite fuel cost column to be called \mosel{FC_f}, where \mosel{f} stands for the actual name of the fuel type. Therefore, the fuel name may not contain special characters and spaces.

To access the database, we have to setup a connection first, so the ODBC data source name is needed again:
\begin{Mosel}{}
if(isExcelDocument(PERIODS_SOURCE)) then
	SQLconnect("DSN=Excel Files;HDR=Yes;DBQ="+expandpath(PERIODS_SOURCE));
	assert(getparam("SQLsuccess"), 'SQLconnect failed!');
else
	SQLconnect(DSN);
	assert(getparam("SQLsuccess"), 'SQLconnect failed!');
end-if
\end{Mosel}

\noindent Note that \mosel{expandpath} expands the path of \mosel{PERIODS_SOURCE}, since the ODBC driver's working directory might be different from ours.

To send an SQL query to the database now, we have to build it up first. The following variables are used for this:
\begin{Mosel}{}
declarations
	periodsIndex: text;						! Index column
	reindexedPeriods: text;					! Period source with converted index
	periodsColumns: list of text;			! Columns to read
	periodsArrays: list of array(K) of real;! Arrays to fill
	periodsQuery: text;						! Assembled SQL query
end-declarations
\end{Mosel}

\noindent The first column to be read is the period index $k$. How to derive this index depends on the indexing of the periods in the database. For our model, we expect the database periods to be indexed by the column \mosel{t} with data type \mosel{TIMESTAMP}. We also assume the database periods to have the same length as the model periods, \mosel{L}; data with different period lengths needs to be interpolated beforehand. 

Given all these informations, we can derive our index $k$ as
\begin{Mosel}{}
periodsIndex := "{fn FLOOR(((t-{ts '"+START+"'})*24+0.1)/"+L+")} + 1";
\end{Mosel}\vskip-1.142em\vskip0pt
\begin{itemize*}
	\item {\ttfamily\smaller \{ts '"+START+"'\}} is the start time \mosel{START} as a timestamp
	\item {\ttfamily\smaller (t-\{ts '"+START+"'\})} is the elapsed time since \mosel{START} in days
	\item {\ttfamily\smaller ((t-\{ts '"+START+"'\})*24+0.1)} is the elapsed time in hours
	\item {\ttfamily\smaller (((t-\{ts '"+START+"'\})*24+0.1)/"+L+")} is the elapsed time in periods $+\frac{L}{10}$
	\item {\ttfamily\smaller \{fn FLOOR(((t-\{ts '"+START+"'\})*24+0.1)/"+L+")\}} is the period index $\in \discrange{0}{T-1}$
	\item {\ttfamily\smaller \{fn FLOOR(((t-\{ts '"+START+"'\})*24+0.1)/"+L+")\} + 1} is our index $k \in \discrange{1}{T}$
\end{itemize*}

\noindent In a fully-fledged database, it would be possible to assign the name \mosel{k} to this converted index, and to use it by this name. In Excel, this assignment has to be done in a subquery. At the same time, we have to specify the source of the periods:
\begin{Mosel}{}
reindexedPeriods := "SELECT "+periodsIndex+" as k, *";
if(isExcelDocument(PERIODS_SOURCE)) then
	reindexedPeriods += " FROM Periods";
else
	reindexedPeriods += " FROM " + PERIODS_SOURCE;
end-if
\end{Mosel}
The subquery \mosel{reindexedPeriods} now represents the periods table, with a converted and renamed index column \verb-k-.

With the index settled, we can select our first data arrays for reading: the demand \mosel{D} and the reserve \mosel{R}.
\begin{Mosel}{}
periodsColumns += ["D", "R"];
periodsArrays += [D, R];
\end{Mosel}
Here, the last statement actually does not copy the \mosel{D} and \mosel{R} arrays. Since they are complex datatypes, Mosel just copies their references to the list. So, when initialising the arrays of this list, we actually initialize our original arrays.

The procedure for the fuel costs is similar, but since the \mosel{FC} array is dynamic, we first have to create the cost subarray for each used fuel type:
\begin{Mosel}{}
forall(f in Fuels) do
	create(FC(f));
	periodsColumns += ["FC_" + f];
	periodsArrays += [FC(f)];
end-do
\end{Mosel}

\noindent Finally, we can assemble the SQL query and specify the needed periods:
\begin{Mosel}{}
periodsQuery := "SELECT k, " + join(periodsColumns, ", ");
periodsQuery += " FROM ("+reindexedPeriods+")";
periodsQuery += " WHERE 1 <= k AND k <= " + T;
\end{Mosel}
The \mosel{join} function joins the elements of \mosel{periodsColumns}, its implementation can be found in Appendix \ref{sec:ExcelFunctions}.

In conclusion, an assembled query for example may equate to
\begin{lstlisting}[language=SQL]
SELECT k, D, R, FC_gas, FC_oil, FC_coal
  FROM (
    SELECT {fn FLOOR(((t-{ts '2009-05-11 00:00:00'})*24+0.1)/1)} + 1 as k, *
      FROM Periods
  )
  WHERE 1 <= k AND k <= 84
\end{lstlisting}
We now initialize our parameters using \mosel{SQLexecute}:
\begin{Mosel}{}
SQLexecute(periodsQuery, periodsArrays);
assert(getparam("SQLsuccess"), 'SQLexecute failed!');
\end{Mosel}

\noindent After the initialization, we have to close the connection to the database server:
\begin{Mosel}{}
SQLdisconnect;
\end{Mosel}

\noindent When working with big Excel spreadsheets, this approach is a good alternative to multiple \mosel{initializations} blocks. Excel reopens the spreadsheet for every \mosel{initializations} block, which can be quite time consuming. In contrast, when using SQL functions, Excel opens the spreadsheet just once at \mosel{SQLconnect}.

\clearpage

\section{Data Output}

Same as with the data input sources, the destination file is given as a parameter:
\begin{Mosel}{}
parameters
	RESULTS_DEST = "Results.xls";	! Destination of the results
end-parameters
\end{Mosel}

\subsection{Electricity price}

Once the cost function has been minimized within the model, we can derive a basic estimate of the electricity price by observing:
\begin{itemize}
  \item The  price is set the highest accepted bid of any unit.
  \item In a market with perfect competition, a unit will bid in at marginal costs.
  \item The marginal cost of a unit is $\varFuel(j) \cdot \fuelCost(\fuelType(j),k) + \varCost(j)$.
\end{itemize}
This can be calculated in Mosel as:
\begin{Mosel}{}
declarations
	price: array(K) of real;
end-declarations
forall(k in K) do
	price(k) := max(j in J |  v(j, k).sol >= 0.5) (FA(j) * FC(F(j), k) + PA(j)); 
end-do
\end{Mosel}

\subsection{Postprocessing}

The maximal production and ramping constraints just impose an upper limit on the maximal possible production. Thus, \mosel{p_max} may actually be lower than the real maximal possible production. We can fix this by deriving the maximal possible production from the production variables and the operational state:
\begin{Mosel}{}
declarations
	exact_p_max: array(J, K) of real;
end-declarations
forall(j in J, k in K) do
	exact_p_max(j, k) := P_max(j).sol * v(j, k).sol;
	if(k > 1) then
		exact_p_max(j, k) := minlist(exact_p_max(j, k),
			p(j, k-1).sol + RU(j) * v(j, k-1).sol
			+ SU(j) * (1 - v(j, k-1).sol) + P_max(j) * (1 - v(j, k).sol));
	end-if
	if(k < T) then
		exact_p_max(j, k) := minlist(exact_p_max(j, k),
			P_max(j) * v(j, k+1).sol + SD(j) * (v(j, k).sol - v(j, k+1).sol));	                       
	end-if
end-do
\end{Mosel}

\subsection{Data Output to Excel}

We want to output the optimizer variables for
\begin{enumerate}
  \item the operational state $\onOff(j,k)$, 
  \item the production variables $\prod(j,k), \maxPossProd(j,k)$,
  \item the price $\price(j)$ and
  \item the costs $\prodCost(j,k), \StartupCost(j,k), \ShutdownCost(j,k)$.
\end{enumerate}

\noindent We will output each variable to its own sheet (inside the same file). This allows us to
\begin{itemize}
  \item change the sizes of $J$ and $K$ and
  \item add new variables to the model
\end{itemize}
without having to change the position of the existing variables. Since there are usually more periods than units, and Excel supports more rows than columns\footnote{Excel 2007 or higher: $2^{20}$ rows vs. $2^{14}$ columns; prior to Excel 2007: $2^{16}$ rows vs. $2^8$ columns or less}, we associate the periods with rows and the units with columns. Together with the actual data, we have to output the index sets \mosel{J} and \mosel{K}. Since the period indices are relative to \mosel{START} and \mosel{L}, we will also output the timestamp of each period.

Now, it would be convenient to implement the output in a separate function, which is to be called for each variable. Unfortunately, one \mosel{initializations} block per function call would be needed, reopening the Excel spreadsheet each time. Depending on the size of the spreadsheet, this may be quite time-consuming. Therefore we use a single \mosel{initializations} block, but still perform the preparation of each variable in a separate function.

We need two functions: one for variables with index set $K$, and one for variables index sets $J$ and $K$. They both should return a two-dimensional \mosel{array of text}, which subsequently will be written to an Excel sheet:
\begin{Mosel}{}
declarations
	PeriodsSheet = dynamic array(rows: range, singleCol: range) of text;
	UnitsPeriodsSheet = dynamic array(rows, cols: range) of text;
end-declarations

function createSheet(data: array(K) of real, name: string): PeriodsSheet
	...
end-function

function createSheet(data: array(J, K) of real): UnitsPeriodsSheet
	...
end-function
\end{Mosel}

\clearpage\noindent In the body of the \mosel{createSheet} function for variables with index set $K$, we need to:
\begin{enumerate}
	\item Fill the first two columns with the index and the starting time of each period:
	\begin{Mosel}{}
forall(k in K) do
	returned(k + 1, 1) := text(k);
	returned(k + 1, 2) := text(datetime(START) + (k-1) * L * 3600);
end-do
	\end{Mosel}
	\item Fill the first row with the given variable name:
	\begin{Mosel}{}
returned(1, 3) := name;
	\end{Mosel}
	\item Fill the space between period indices and titles with the actual data:
	\begin{Mosel}{}
forall(k in K) do
	returned(k + 1, 3) := text(data(k));
end-do
	\end{Mosel}
\end{enumerate}

\noindent In the function for variables with index sets $J$ and $K$, the time series for all units $J$ should be written to adjacent columns. Since the index set $J$ is not necessarily of the form $\{1, \dots, |J|\}$, we first construct a mapping from $J$ to $\{1, \ldots, |J|\}$:
\begin{Mosel}{}
declarations
	mapping: dynamic array(J) of integer;
end-declarations
forall(j in J) do
	mapping(j) := getsize(mapping) + 1;
end-do	
\end{Mosel}

\noindent The remainder of the function is similiar to the previous function for index set $K$. We need to:
\begin{enumerate}
	\item Fill the first two columns with the index and the starting time of each period:\begin{Mosel}{}
forall(k in K) do
	returned(k + 1, 1) := text(k);
	returned(k + 1, 2) := text(datetime(START) + (k-1) * L * 3600);
end-do
\end{Mosel}
	\item Fill the first row with the indices of the units,
	\begin{Mosel}{}
forall(j in J) do
	returned(1, mapping(j)) := text(j);
end-do
	\end{Mosel}
	\item Fill the space between period indices and titles with the actual data:
	\begin{Mosel}{}
forall(j in J, k in K) do
	returned(k + 1, mapping(j)) := text(data(j, k));
end-do
\end{Mosel}
\end{enumerate}

\noindent These functions only take \mosel{real}-valued time series; For convenience, we implement them for \mosel{mpvar}-valued time series as well:
\begin{Mosel}{}
function createSheet(data: array(K) of mpvar, name: string): PeriodsSheet
	returned := createSheet(array(k in K) data(k).sol, name);
end-function

function createSheet(data: array(J, K) of mpvar): UnitsPeriodsSheet
	returned := createSheet(array(j in J, k in K) data(j, k).sol);
end-function
\end{Mosel}

\noindent To perform the output to Excel, we still need to specify the region to which each variable needs to be written. For this, we will use the \mosel{cellRange} function, implemented in our ExcelFunctions package (Appendix \ref{sec:ExcelFunctions}), which generates Excel ranges like \verb-[Sheet$A1:Q10]-.
\begin{Mosel}{}
initializations to "mmodbc.excel:noindex;" + RESULTS_DEST
	evaluation of createSheet(v)
		as cellRange("v", 1, 1, getsize(J)+2, getsize(K)+1);
	evaluation of createSheet(p)
		as cellRange("p", 1, 1, getsize(J)+2, getsize(K)+1);
	evaluation of createSheet(exact_p_max)
		as cellRange("p_max", 1, 1, getsize(J)+2, getsize(K)+1);
	evaluation of createSheet(price, "Price")
		as cellRange("price", 1, 1, 3, getsize(K)+1);
	evaluation of createSheet(s)
		as cellRange("s", 1, 1, getsize(J)+2, getsize(K)+1);
	evaluation of createSheet(c)
		as cellRange("c", 1, 1, getsize(J)+2, getsize(K)+1);
	evaluation of createSheet(cp)
		as cellRange("cp", 1, 1, getsize(J)+2, getsize(K)+1);
	evaluation of createSheet(cu)
		as cellRange("cu", 1, 1, getsize(J)+2, getsize(K)+1);
	evaluation of createSheet(cd)
		as cellRange("cd", 1, 1, getsize(J)+2, getsize(K)+1);
end-initializations
\end{Mosel}

\clearpage

\section{Remarks}

This whitepaper presents a detailled model for the Unit Commitment problem. We have given its practical implementation and discussed how to reduce the number of start-up cost constraints and how to analyze the cause of infeasibilities. Finally, we have shown how to manage data input and output, and have given a basic estimate of the electricity price.

One major drawback of the implemented model is the lacking consideration of market power. To remedy this, usually a so-called ``uplift'' function is added to the model's electricity price, which is derived by statistically analyzing the disagreement between the outcome of this model and the real price on historical data.

Of course, it would be better to consider the market power directly as part of the model. Different approaches in achieving this are discussed in \cite{ventosa_electricity_2005}. The two most popular ones are game theoretical approaches:
\begin{itemize}
  \item the Cournot equilibrium and
  \item the Supply Function equilibrium.
\end{itemize}
They differ in the modelling of the strategies of the power generating companies. While in a Cournot model a company may only decide on it's power output, a company's strategy in a Supply Function equilibrium is described by a function mapping the market price to its power supply. The greater freedom in the choice of a strategy, however, comes at the expense of a higher computational effort, which in turn forces the use of less detailled models.

For a comparison of the practical performance of Cournot and Supply Function models on the German power market, we refer to \cite{willems_cournot_2009}.

The second drawback is the neglection of the underlying power grid. The production schedule resulting from our model may not be feasible on a real power grid, i.e. the power grid may not be able to transport the electricity from the producing units to the consumers. While this was typically not a fundamental issue so far, it becomes more and more important due to the increasing energy production from renewable, more volatile resources causing bottlenecks in the power exchange beetween different regions and thus necessitate an explicit modelling of the power grid.

Also, most of today's power markets are connected to one or more neighboring markets through interconnectors. The interconnectors are used to transport energy from a market with higher price to a market with lower price, thus diminishing the price difference. The effective price change is typically small due to the small capacity of the interconnectors, but may not be negligible.

\clearpage

\bibliography{Whitepaper}

\clearpage


\appendix

\section{Excel Functions}
\label{sec:ExcelFunctions}

The Excel functions are implemented as part of the ExcelFunctions package. Therefore, the functions to be used in other packages or the main model have to marked as \mosel{public}.

\subsection{Function isExcelDocument}

The function \mosel{isExcelDocument} checks if the given data source is an Excel
file. This is done by checking the file extension:
\begin{Mosel}{}
public function isExcelDocument(source: string): boolean
	declarations
		S = getsize(source);
	end-declarations
	returned := strtolower(copytext(text(source), S-4+1, S)) = '.xls' or
	            strtolower(copytext(text(source), S-5+1, S)) = '.xlsx';
end-function
\end{Mosel}

\noindent Unfortunately, we also have to implement the \mosel{strtolower} function for ourselfs, which converts a string to lower case:
\begin{Mosel}{}
public function strtolower(str: text): text
	returned := str;
	forall(i in 1..getsize(returned)) do
		if(65(! = A !) <= getchar(returned, i) and
		   getchar(returned, i) <= 90(! = Z !)) then
			setchar(returned, i, getchar(returned, i) + (97-65));
		end-if
	end-do
end-function
\end{Mosel}

\subsection{Referencing cells in Excel}

Excel supports two notations for referencing cells,
\begin{itemize}
  \item the A1 notation (ex: C4, R7, Z80) and
  \item the R1C1 notation (ex: R4C3, R7C18, R80C26).
\end{itemize}
While the R1C1 notation is easier to use programmatically, the A1 notation may be more familiar to end users. Therefore, we offer the boolean \mosel{useA1notation} which switches from R1C1 to A1 notation:
\begin{Mosel}{}
public declarations
	useA1notation: boolean;	! If true, the A1 notation is used
end-declarations
\end{Mosel}

Depending on this setting, the coordinates of a cell can be determined with a call to \mosel{cellCoords}:
\begin{Mosel}{}
public function cellCoords(column: integer, row: integer): string
	if(useA1notation) then
		returned := columnToA1(column) + row;
	else
		returned := "R" + row + "C" + column;
	end-if
end-function
\end{Mosel}
($\rightarrow$ \mosel{columnToA1} is described in the next subsection)

\noindent Cell ranges are described by the coordinates of their top-left and bottom-right cells, divided by a colon (ex. C4:R7, R4C3:R7C18). They can be determined with a call to \mosel{cellRangeCoords}:
\begin{Mosel}{}
public function cellRangeCoords(column: integer, row: integer,
                                width: integer, height: integer): string
	returned := cellCoords(column, row) + ":"
	          + cellCoords(column + width - 1, row + height - 1);
end-function
\end{Mosel}

\noindent When referencing Excel cells through ODBC, it is usual to include the cell's sheet name too. This is done by prepending the sheet name to the coordinates, divided by a dollar sign, and by enclosing the reference in square brackets. The two following functions generate such references to a cell or respectively to a cell range:
\begin{Mosel}{}
public function cell(sheet: string, column: integer, row: integer): string
	returned := "[" + sheet + "$" + cellCoords(column, row) + "]";
end-function
public function cellRange(sheet: string, column: integer, row: integer,
                          width: integer, height: integer): string
	returned := "["+sheet+"$"+cellRangeCoords(column, row, width, height) + "]";
end-function
\end{Mosel}

\subsubsection{Column index conversions in the A1 notation}

The column numbering scheme in the A1 notation is a bit unusual and needs special conversion functions.

\noindent Conversion from integer to A1 column index (1 $\rightarrow$ A, 2 $\rightarrow$ B, ...)
\begin{Mosel}{}
public function columnToA1(index: integer): string
	declarations
		c: text;		! Column name as text (needs mmsystem)
	end-declarations
	while(index > 0) do
		c := " " + c;
		setchar(c, 1, 65 (! = A !) + ( (index-1) mod 26));
		index := (index-1) div 26;
	end-do
	returned := string(c);
end-function
\end{Mosel}

\noindent Conversion from A1 column index to integer (A $\rightarrow$ 1, B $\rightarrow$ 2, ...)
\begin{Mosel}{}
public function a1ToColumn(column: string): integer
	forall(i in 1..getsize(column)) do
		returned := returned * 26 + getchar(column, i) - 65 (! = A !) + 1;
	end-do
end-function
\end{Mosel}

\subsection{Function Join}

Joins the specified separator string between each element of the specified set of strings, yielding a single joined string. Example:
\begin{Mosel}{}
join(["A", "D", "B", "C"], ", ") -> "A, D, B, C"
join(["A"], ", ") -> "A"
join({"A", "D", "B", "C"}, ", ") -> "A, B, C, D"
\end{Mosel}
Implementations for lists and sets of string:
\begin{Mosel}{}
public function join(strings: list of text, separator: string): text
	forall(s in strings) do
		returned := returned + s + separator;
	end-do
	returned -= separator;
end-function

public function join(strings: set of text, separator: string): text
	forall(s in strings) do
		returned := returned + s + separator;
	end-do
	returned -= separator;
end-function
\end{Mosel}

\end{document}